\documentclass[sn-mathphys,Numbered]{sn-jnl}

\usepackage{graphicx}%
\usepackage{multirow}%
\usepackage{amsmath,amssymb,amsfonts}%
\usepackage{amsthm}%
\usepackage{mathrsfs}%
\usepackage[title]{appendix}%
\usepackage{xcolor}%
\usepackage{cases}
\usepackage{textcomp}%
\usepackage{manyfoot}%
\usepackage{booktabs}%
\usepackage{indentfirst}
\usepackage{algpseudocode}%
\usepackage{listings}%
\usepackage{hyperref}
\usepackage[ruled,vlined]{algorithm2e}
\usepackage{setspace}

\begin{document}

\title[Article Title]{A New Fast Adaptive Linearized Alternating Direction Multiplier Method for Convex Optimization}


\author{\fnm{Boran Wang}}\email{1932210936@qq.com}

\affil{\orgdiv{College of Science}, \orgname{Minzu University of China},  \city{Beijing}, \postcode{100081}, \state{China}}


\abstract{This work proposes a novel adaptive linearized alternating direction multiplier method (LADMM) to convex optimization, which improves the convergence rate of the LADMM-based algorithm by adjusting step-size iteratively.
The innovation of this method is to utilize the information of the current iteration point to adaptively select the appropriate parameters, thus expanding the selection of the subproblem step size and improving the convergence rate of the algorithm while ensuring convergence.
The advantage of this method is that it can improve the convergence rate of the algorithm as much as possible without compromising the convergence. This is very beneficial for the solution of optimization problems because the traditional linearized alternating direction multiplier method has a trade-off in the selection of the regular term coefficients: larger coefficients ensure convergence but tend to lead to small step sizes, while smaller coefficients allow for an increase in the iterative step size but tend to lead to the algorithm's non-convergence. This balance can be better handled by adaptively selecting the parameters, thus improving the efficiency of the algorithm.}

\keywords{variational inequality, Linearized ADMM algorithms,self-adaptation, global convergence.}



\maketitle

\section{Introduction}\label{sec1}
 
The convex minimization problem under consideration involves linear constraints and an objective function expressed as the sum of two functions that do not involve coupled variables.
\begin{equation}\label{eqn:1} 
\min\left \{\theta _{1} \left ( x \right ) +  \theta _{2} \left ( y \right ) | Ax+By=b,x\in X,y\in Y  \right \} ,
\end{equation}
where $\theta _{1} \left (x \right):R^{n_{1} } \to R$ and $\theta _{2} \left (x \right):R^{n_{2} } \to R$ are convex functions (but not necessarily smooth), $A\in R^{m\times n_{1} } ,B\in R^{m\times n_{2} } ,b\in R^{m},X\subseteq R^{n_{1} }$ and$Y\subseteq R^{n_{2} }$ are closed convex sets.Throughout this study, we assume that matrices A and B are of full column rank. The problem \eqref{eqn:1} has been utilized in several fields, especially image processing \cite{1,2}, statistical learning \cite{3}, and communication networking \cite{4,5}.

The augmented Lagrangian function of Problem \eqref{eqn:1} is
\begin{equation}\label{eqn:2}
 L_{\beta }\left ( x,y,\lambda  \right )=\theta _{1}\left ( x \right) +\theta _{2}\left ( y \right )-\lambda ^{T}\left ( Ax+By-b \right )+\frac{\beta }{2}\left \| Ax+By-b \right \|_{2}^{2}.   
\end{equation}
where $\lambda \in R^{m}$ represents the Lagrange multiplier, and $ \beta > 0$ represents the penalty parameter.

The augmented Lagrangian method (ALM) created in \cite{6,7} is a classic and essential algorithm for addressing \eqref{eqn:1}. It plays an important role in both theoretical research and algorithmic design for diverse convex programming problems. For a better understanding of the ALM, see \cite{8,9,10,11,12,13,14}. The alternate direction method of multipliers (ADMM), a split version of the augmented Lagrangian approach first proposed in \cite{15,16}, is useful for solving (1). For a better understanding of the ALM, see \cite{17,18,19}. The ADMM decomposes the augmented Lagrangian independently for the variables $x$ and $y$, allowing the subproblems to be addressed separately. In every iteration, the  ADMM proceeds as follows in each iteration:
\begin{equation}\label{3}
\left\{
\begin{aligned}
x^{k+1} &=arg\min_{x}\left \{ L_{\beta } \left ( x,y^{k},\lambda ^{k}   \right )|x\in X  \right \}
 \\y^{k+1} &=arg\min_{y}\left \{L_{\beta } \left ( x^{k+1} ,y,\lambda ^{k}   \right )| y\in Y \right \}
 \\\lambda ^{k+1}& =\lambda ^{k} -\beta \left ( Ax^{k+1} +By^{k+1} -b \right ) 
\end{aligned}
\right.
\end{equation}

The ADMM has been studied extensively in the literature, with a focus on solving subproblems more efficiently for various scenarios. Special properties or structures of the functions $\theta _{1}$ and $\theta _{2}$, or the coefficient matrices $A$ and $B$, may be helpful in designing application-tailored algorithms based on the prototype ADMM scheme \eqref{3}. In theory, convergence should still be ensured. This concept emphasizes the significance of applying the ADMM to specific applications in different fields. For clarification, consider the $y -$ subproblem in \eqref{3}, which may be expressed as
\begin{equation}\label{4}
 y^{k+1}=arg\min_{y}\left \{\theta_{2}\left ( y \right ) +  \frac{\beta }{2} \left \| By+\left ( Ax^{k+1}-b-\frac{1}{\beta } \lambda ^{k}   \right )  \right \|^{2}  | y\in Y \right \}   
\end{equation}

The above illustrates that the function $g(y),$ the matrix $B$, and the set $Y$ are necessary for solving the $y-$subproblem. When the function, matrix, and set are in the generic form, we can solve \eqref{4} using standard iteration\cite{20,21}. The functions, matrices, and sets $Y$ in the $y-$ subproblems are unique in some real-world scenarios, necessitating more efficient solutions. \\
After linearizing the quadratic term $\left \| By+ \left ( Ax^{k+1} -b-\frac{1}{\beta }\lambda ^{k} \right ) \right \|^{2}$ in \eqref{4}, it is obtained:
\begin{equation}
y^{k+1}=arg\min_{y}\left \{ \theta_{2}\left ( y \right )+  \frac{r}{2} \left \| y-\left ( y^{k}+  \frac{1}{r}q^{k}    \right )  \right \| ^{2} |y\in  Y    \right \}.    
\end{equation}
where
\begin{equation}
q^{k}=B^{T}\left ( \lambda ^{k}-\beta \left ( Ax^{k+1}+By^{k}-b   \right )   \right ).    
\end{equation}
where $ r> 0 $ is a constant, and thus we obtain the linearized alternating direction multiplier method:
\begin{equation}
\left\{
\begin{aligned}
x^{k+1} &=arg\min_{x}\left \{ L_{\beta } \left ( x,y^{k},\lambda ^{k}   \right )|x\in X  \right \}
 \\y^{k+1} &=arg\min_{y}\left \{L_{\beta } \left ( x^{k+1} ,y,\lambda ^{k}   \right )+\frac{r}{2}\left \| y-\left (  y^{k}+\frac{1}{r}q^{k}   \right )  \right \|^{2}   | y\in Y \right \}
 \\\lambda ^{k+1}& =\lambda ^{k} -\beta \left ( Ax^{k+1} +By^{k+1} -b \right ) 
\end{aligned}
\right.
\end{equation}
Linearized alternating direction multiplier method is widely used in areas; we refer to , e.g.\cite{22,23,24,25} . 

A generalized alternating direction method of multipliers\cite{26} can be written as follows:
\begin{equation}\label{8}
\left\{
\begin{aligned}
x^{k+1} &=arg\min_{x}\left \{ L_{\beta } \left ( x,y^{k},\lambda ^{k}   \right )|x\in X  \right \}
 \\y^{k+1} &=arg\min_{y}\left \{L_{\beta } \left ( x^{k+1} ,y,\lambda ^{k}   \right )+\frac{1}{2}\left \| y-y^{k}    \right \|_{D} ^{2}   | y\in Y \right \}
 \\\lambda ^{k+1}& =\lambda ^{k} -\beta \left ( Ax^{k+1} +By^{k+1} -b \right ) 
\end{aligned}
\right.
\end{equation}
where $D\in R^{n_{2} \times n_{2} }$ is positive definite, and thus, the linearized alternating direction multiplier method can be viewed as a special case of the general  alternating direction multiplier method\cite{26}, where the positive definite proximal term.
\begin{equation}
D=rI_{n_{2} }-\beta B^{T}B\quad and \quad r> \beta \left \| B^{T}B  \right \|    
\end{equation}
Since in many practical applications, only one subproblem  needs to be linearized, in this paper, we only need to consider the case of linearizing the $y-$ subproblem in  \eqref{8}.

Literature\cite{27} proposes an indefinite proximity alternating direction multiplier method with the following iteration format:
\begin{equation}
\left\{
\begin{aligned}
x^{k+1} &=arg\min_{x}\left \{ L_{\beta } \left ( x,y^{k},\lambda ^{k}   \right )|x\in X  \right \}
 \\y^{k+1} &=arg\min_{y}\left \{L_{\beta } \left ( x^{k+1} ,y,\lambda ^{k}   \right )+\frac{1}{2}\left \| y-y^{k}    \right \|_{D_{0} } ^{2}   | y\in Y \right \}
 \\\lambda ^{k+1}& =\lambda ^{k} -\beta \left ( Ax^{k+1} +By^{k+1} -b \right ) 
\end{aligned}
\right.
\end{equation}
where  $$D_{0} =\tau rI_{n_{2} }-\beta B^{T}B,\quad0.8\le\tau < 1\quad  and \quad r> \beta \left \| B^{T}B  \right \|.$$

He et al. proposed an indefinite proximity linearized ADMM in the literature \cite{28} with the iterative steps:
\begin{equation}
\left\{
\begin{aligned}
x^{k+1} &=arg\min_{x}\left \{ L_{\beta } \left ( x,y^{k},\lambda ^{k}   \right )|x\in X  \right \}
 \\y^{k+1} &=arg\min_{y}\left \{L_{\beta } \left ( x^{k+1} ,y,\lambda ^{k}   \right )+\frac{1}{2}\left \| y-y^{k}  \right \|_{D}^{2}   | y\in Y \right \}
 \\\lambda ^{k+1}& =\lambda ^{k} -\beta \left ( Ax^{k+1} +By^{k+1} -b \right ) 
\end{aligned}
\right.
\end{equation}
where$$D=\tau r I-\beta B^{T}B, \quad0.75< \tau < 1\quad and \quad\rho > \beta \left \| B^{T}B  \right \| $$
It is easy to see that the matrix D is not necessarily semipositive definite. In the paper the authors prove that 0.75 is an optimal lower bound for $\tau$.

Based on the above discussion, the above linearized alternating direction multiplier methods are fixed in the selection of the parameter $\tau$. However, the fixed $\tau$ may be computationally conservative. This paper adopts an adaptive technique to adaptively select the values of the regular term coefficients $\delta _{k}$ and proposes an adaptive linearized alternating direction multiplier method. The regular term coefficients $\delta _{k}$ are determined automatically based on the knowledge of the current iteration point, resulting in a more efficient process. The approach has a broader application and improved solution efficiency. This work provides a comprehensive convergence analysis for the proposed method using variational inequalities and optimization theory, as well as numerical tests to verify the approach's efficacy.

The remainder of the paper is arranged as follows. Section II provides the necessary background information for the theoretical analysis; Section III suggests a novel fast adaptive linearized alternating direction multiplier approach, as well as a convergence analysis; Section IV presents the numerical experiments; and Section VII makes conclusions.

\section{Preliminaries}\label{sec2}
\noindent We summarize some preliminary findings and present some straightforward results in this part that will be used to the convergence analysis.

\noindent In this article, the symbol $\left \| \cdot \right \|$ denotes the two-norm $\left \| \cdot \right \|_{2} .\left \| x \right \|_{D}^{2}: = x^{T}Dx$ is the matrix norm, where $D\in R^{n\times n} $ is a symmetric positive definite matrix and the vector $x\in R^{n}$. When $D $ is not a positive definite matrix, we still use the above notation.
\subsection{characterization of variational inequalities}\label{subsec2}

Let the Lagrangian function of \eqref{eqn:1} defined on $X\times Y\times R^{m}$ be :
\begin{equation}\label{eqn:112}
 L\left ( x ,y,\lambda   \right ) =\theta _{1} \left ( x \right )+\theta _{2} \left ( y \right ) -\lambda ^{T} \left ( Ax+By-b \right ).
 \end{equation} 
A pair $\left ( \left ( x^{\ast },y^{\ast }   \right ),\lambda ^{\ast }   \right )$ is called a saddle point of the Lagrangian function \eqref{eqn:112} if it satisfies the inequalities
$$  L_{\lambda \in R^{m} } \left ( x^{\ast } ,y^{\ast } ,\lambda  \right )\le L\left ( x^{\ast } ,y^{\ast } ,\lambda ^{\ast }  \right )\le L_{x\in X,y\in Y} \left ( x ,y ,\lambda^{\ast }   \right ) .$$
They can be reformulated as variational inequalities:
\begin{equation}\label{eqn:113}
\begin{cases}
x^{\ast } \in  X,\qquad \theta _{1} \left ( x \right )- \theta _{1} \left ( x^{\ast }  \right )+\left ( x-x^{\ast }  \right )\left ( -A^{T} \lambda ^{\ast }  \right )\ge 0,\quad\forall x\in X . 
\\
y^{\ast } \in  Y ,\qquad \theta _{2} \left ( y \right )-  \theta _{2} \left ( y^{\ast }  \right )+\left ( y-y^{\ast }  \right )\left ( -A^{T} \lambda ^{\ast }  \right )\ge 0,\quad\forall y\in Y .
\\\lambda ^{\ast } \in R^{m} ,\qquad\qquad\left ( \lambda -\lambda ^{\ast }  \right )^{T}\left ( Ax^{\ast } +  By^{\ast } -b \right ) \ge 0,\quad\forall \lambda \in R^{m} .  
\end{cases}
\end{equation}
or in its condensed form:
\begin{subequations}\label{eqn:114}
\begin{align}
    & VI(F,\theta ,\Omega ):\quad w^{\ast } \in \Omega ,\quad \theta \left ( u \right ) -\theta \left ( u^{\ast }  \right )+  \left ( w-w^{\ast }  \right ) ^{T}F\left ( w^{\ast }  \right ) \ge 0,\quad \forall w\in \Omega . \label{eqn:114a} \\
    &  \begin{aligned}
        \intertext{where} \\
         & \theta \left ( u \right ) =\theta _{1} \left ( x \right ) + \theta _{2} \left ( y \right ) , 
          \quad u=\begin{pmatrix}x \\y \end{pmatrix}, 
          \quad w=\begin{pmatrix}x \\y \\\lambda \end{pmatrix}, 
          \quad F\left ( w \right ) =\begin{pmatrix}-A^{T}\lambda \\ -B^{T}\lambda \\ Ax+By-b \end{pmatrix}
    \end{aligned} \label{eqn:114b}
\end{align}
\end{subequations}

and $$\Omega =X\times Y\times R^{m} .$$
We denote by $\Omega ^{\ast } $ the solution set of \eqref{eqn:114}. Note that the operator $F$ in \eqref{eqn:114b} is affine with a skew-symmetric matrix. Thus we have 
\begin{equation}
    \left ( w-\bar{w}  \right )^{T}  \left ( F\left ( w \right )-F\left ( \bar{w}  \right )   \right ) = 0,\quad\forall w,\bar{w}.
\end{equation}

The following lemma will be frequently used later. Its proof can be found in the literature, see, e.g., Theorem 3.71 in \cite{29}.\\
\textbf{Lemma 2.1}. Let $X\subset R^{n}$ be a closed convex set, $\theta \left ( x \right )$ 
 and $ f \left ( x \right ) $ be convex functions. Assume that $f$ is differentiable and the solution set of the problem $min\left \{ \theta \left ( x \right )+f\left ( x \right )   \right \}$ is  nonempty. Then we have
\begin{subequations}
\begin{align}
     &x^{\ast } \in \text{argmin}\left \{ \theta \left ( x \right )+f\left ( x \right ) \mid x\in X  \right \} , \\
\intertext{if and only if}
x^{\ast } \in X, \quad & \theta \left ( x \right ) -\theta \left ( x^{\ast }  \right ) +\left ( x-x^{\ast }  \right )^{T}\nabla f\left ( x^{\ast }  \right )\ge 0, \quad \forall x\in X.
\end{align}
\end{subequations}

\subsection{some notation}\label{subsec2}

We define various variables and matrices to make analysis easier. Let

\begin{equation}\label{eqn:117}
 Q_{k+1}=\begin{pmatrix}
  \beta   \delta _{k}I_{n_{2} }  & O\\
  -B&\frac{1}{\beta }I_{m }  
\end{pmatrix}  \qquad and \qquad M=\begin{pmatrix}
  I_{n_{2} } &O \\
  -\beta B&I_{m }
\end{pmatrix}.
\end{equation}
\begin{equation}\label{eqn:118}
H_{k+1}=\begin{pmatrix}
  \beta   \delta _{k}I_{n_{2} }  & O\\
  O&\frac{1}{\beta } I_{m } 
\end{pmatrix} 
\end{equation}

\begin{equation}\label{eqn:119}
\tilde{w}^{k}=\begin{pmatrix}\tilde{x}^{k} 
 \\\tilde{y}^{k}
 \\\tilde{\lambda}^{k}

\end{pmatrix} =\begin{pmatrix}x^{k+1} 
 \\y^{k+1}
 \\\lambda  ^{k}  -\beta \left ( Ax^{k+1} +By ^{k} -b\right ) 
\end{pmatrix}
\end{equation}
\begin{equation}\label{eqn:120}
\tilde{v}^{k}=\begin{pmatrix}
 \tilde{y}^{k}
 \\\tilde{\lambda}^{k}

\end{pmatrix} =\begin{pmatrix} 
 y^{k+1}
 \\\lambda  ^{k}  -\beta \left ( Ax^{k+1} +By ^{k} -b\right ) 
\end{pmatrix}    
\end{equation}
\textbf{Lemma 2.2}. The $Q_{k+1},H_{k+1}$ and $M$ defined in \eqref{eqn:117} and \eqref{eqn:118} satisfy:
\begin{equation}\label{eqn:121}
 Q_{k+1}=H_{k+1}M\quad and\quad  H_{k+1}\succeq0,   
\end{equation}
Proof: \eqref{eqn:121} clearly hold. \\
\textbf{Lemma 2.3}. (Robbins-Siegmund Theorem\cite{30}) $a^{k}, b^{k}, c^{k} $ and $d^{k}$ are non-negative sequences and there are :
\begin{equation}\label{eqn:122}
a^{k+1}\le \left ( 1+b^{k}  \right ) a^{k}+c^{k}-d^{k} , \forall k=0,1,2\dots    
\end{equation}
if  $\sum_{k=0}^{+  \infty }b^{k} < +  \infty$ and $\sum_{k=0}^{+  \infty }c^{k} < +  \infty,$ so $\lim_{k \to \infty} a^{k}$ exists and is bounded, while there are $\sum_{k=0}^{+  \infty }d^{k} < +  \infty.$\\
\section{Algorithm and convergence analysis}\label{sec3}
\subsection{New algorithms}\label{subsec2}
Suppose that $\theta _{1} :R^{n_{1} }\to R\cup \left \{ +\infty \right \}$ and $\theta _{2}:R^{n_{2}}\to R\cup \left \{ +\infty \right \}$ are proper, closed, convex functions. \\
Let 
\begin{equation}\label{eqn:123}
H\left ( y \right ) = \frac{1}{2} \left \| Ax^{k+1}+By-b   \right \|^{2}.    
\end{equation}
The expansion of $H$ in a Taylor series around $y^{k}$ can be expressed
\begin{equation}\label{eqn:124}
\begin{aligned}
\frac{1}{2} \left \| Ax^{k+1}+By-b   \right \|^{2}&\approx \frac{1}{2}  \left \| Ax^{k+1}+By^{k} -b   \right \|^{2}
\\&+  \left \langle B^{T}\left ( Ax^{k+1}+By^{k} -b \right ),y-y^{k}    \right \rangle+\frac{\mu _{k} }{2}\left \| y-y^{k}  \right \|^{2}.      
\end{aligned}
\end{equation}
where
\begin{equation}\label{eqn:125}
\begin{aligned}
\mu _{k} = \left \| \bigtriangledown H\left ( y^{k}  \right ) -\bigtriangledown H\left ( y^{k-1}  \right )-\mu \left ( y^{k}-y^{k-1}   \right )  \right \|^{2}  = \frac{\left \| B\left( y^{k}-y^{k-1} \right) \right \|^{2}}{\left \| \left( y^{k}-y^{k-1} \right) \right \|^{2}}.      
\end{aligned}
\end{equation}

\noindent The main iterative steps of the algorithm are:
\begin{spacing}{0.9}
\noindent\makebox[\textwidth]{\rule[0.1ex]{\textwidth}{0.01pt}}\\
Algorithm 1. A New Fast Adaptive Linearized Alternating Direction Multiplier Method for Convex Optimization.\\[-1ex]
\makebox[\textwidth]{\rule[0.1ex]{\textwidth}{0.01pt}}
\end{spacing}
\noindent \textbf{Step 0.} Choose parameters $\tau > 1 $, $ \eta > 1$, $\beta > 0$, $\delta _{min} > 0$, $\varepsilon \in \left ( 0,\frac{1}{2}  \right )$ and set $\delta _{0} =\delta _{-1} =0.75\left \| B^{T}B  \right \|$ . Select an initial point $ w^{0}=\left ( x^{0},y^{0},\lambda ^{0}    \right ) \in \Omega$, where $\Omega  = X\times Y\times R^{m}$. set up $k=0. $  \\
\textbf{Step 1.} Calculate: $w^{k+1}=\left ( x^{k+1},y^{k+1},\lambda ^{k+1}    \right ) \in X\times Y\times R^{m}.$  
\begin{subequations}\label{eqn:126}
\begin{numcases}{\mbox{}}
x ^{k+1}=arg\min_{x}\left \{\theta _{1} \left ( x\right ) -\left (\lambda ^{k}   \right )^{T}Ax  +  \frac{\beta }{2} \left \| Ax+By^{k} -b  \right \|^{2}  \right \}  \label{eqn:126:1}\\
y ^{k+1}=arg\min_{y}\left \{\theta _{2} \left ( y\right )-\left (\lambda ^{k}   \right )^{T}By+\beta \left \langle B^{T}\left ( Ax^{k+1}+By^{k}-b   \right ),y-y^{k}    \right \rangle\right.\notag\\\left.\qquad\qquad\qquad\qquad\qquad +  \frac{  \delta _{k}\beta }{2}\left \| y-y^{k}  \right \| ^{2}   \right \} \label{eqn:126:2}\\
\lambda ^{k+1}    =\lambda ^{k}-\beta \left (  Ax ^{k+1}+By ^{k+1}-b  \right ) \label{eqn:126:3}\end{numcases}
\end{subequations}
\textbf{Step 2.} If any one of the following conditions holds:
\begin{align}\label{eqn:127}
Condition 1.\quad & \delta_{k} \left \| y^{k+1}-y^{k} \right \|^{2} > \frac{1}{2\varepsilon} \left \| B\left( y^{k+1}-y^{k} \right) \right \|^{2} . \\
Condition 2.\quad & y^{k+1}=y^{k}. \nonumber
\end{align}
then go to step 3. otherwise, set up $\delta _{k}  =\tau \ast\delta _{k}\left ( \tau > 1 \right )  ,$ go back to step 1.\\
\textbf{Step 3.} If $\delta _{k}>  \delta _{k-1}$ , then   $ \delta _{min} $ is replaced by $   \eta \delta _{min} $ . \\
\textbf{Step 4.} Calculate $\delta _{k+1}$: \begin{equation}\label{eqn:128}
\delta _{k+1} =\max\left \{  h _{k+1},\min\left ( \delta _{min},\left \| B^{T}B  \right \|   \right )    \right \}  . 
\end{equation}
where
\begin{equation*}
 h _{k+1} = \begin{cases} \delta_{k}, & \text{if } y^{k+1} = y^{k}. \\ \frac{\left \| B\left ( y^{k+1} -y^{k}  \right )  \right \|^{2}  }{\left \|  y^{k+1} -y^{k}    \right \|^{2}} , & \text{if } y^{k+1} \neq y^{k}. \end{cases} 
\end{equation*}
\textbf{Step 5.} If the stopping criterion is satisfied, terminate the algorithm.\\
otherwise, set $k+1=k,$ and go to step 1.

\begin{spacing}{0.9}
\noindent\makebox[\textwidth]{\rule[0.1ex]{\textwidth}{0.01pt}}\\
\end{spacing}
\subsection{Global convergence analysis}\label{subsec2}
\noindent For the problem \eqref{eqn:1}, we establish the global convergence analysis of Algorithm 1 in this section. First, let's establish four essential lemmas.

We further refer to the following notation in order to discern the crucial variables:
\begin{equation}\label{eqn:129}
    v=\begin{pmatrix}y
 \\\lambda 
\end{pmatrix},\quad V=Y\times R^{m} ,\quad and \quad V^{\ast } =\left \{ \left ( y^{\ast },\lambda ^{\ast }   \right )|\left ( x^{\ast } ,y^{\ast },\lambda ^{\ast }   \right )\in \Omega ^{\ast }    \right \}
\end{equation}
\textbf{Lemma 3.1}. Let  $ \left \{ w^{k} \right \} $ be the sequence generated by 
 the Algorithm 1 for the problem \eqref{eqn:1} and $ \tilde{w} ^{k}  $ be defined by\eqref{eqn:119} . Then we have $ \tilde{w} ^{k} \in \Omega $ and
\begin{equation}\label{eqn:130}
\theta \left ( u \right ) -\theta \left ( \tilde{u} ^{k}  \right ) +\left ( w-\tilde{w}^{k}   \right )^{T} F\left ( \tilde{w}^{k} \right )  \ge \left ( v-\tilde{v}^{k}   \right )^{T}  Q_{k+1} \left (  v^{k}-\tilde{v}^{k}   \right ),\quad\forall w\in \Omega . 
\end{equation}
where the matrix $Q_{k+1}$ is defined  in \eqref{eqn:117}.\\
Proof: First, from the optimality conditions of the subproblems \eqref{eqn:126:1} and \eqref{eqn:126:2} , we respectively have 
\begin{equation}\label{eqn:131}
\begin{split}
&x^{k+1}\in X ,\quad  \theta _{1} \left ( x \right )  -\theta _{1} \left ( x ^{k+1}  \right ) +\left ( x-x ^{k+1}  \right )^{T}\left \{ -A^{T}\lambda  ^{k}  +\beta A^{T}\left ( Ax^{k+1} +By^{k} -b   \right ) \right \} 
\\&\quad\ge 0,\quad \forall x\in X, 
\end{split}
\end{equation}
and
\begin{equation}\label{eqn:132}
\begin{split}
y^{k+1}\in Y ,\quad \theta _{2}\left ( y \right ) -\theta _{2}\left ( y^{k+1}  \right )&+\left ( y-y^{k+1}   \right )^{T}\left \{-B^{T}\lambda ^{k}+  \beta B^{T}\left ( Ax^{k+1}+  By^{k} -b   \right )\right.\\& \left.+  \beta   \delta _{k}\left ( y^{k+1}-y^{k}    \right )  \right \}\ge 0,\quad\forall y\in Y.
\end{split} 
\end{equation}
Based on \eqref{eqn:119} and \eqref{eqn:120}, we get
\begin{equation}\label{eqn:133}
\tilde{x}^{k} =x^{k+1},\quad \tilde{y}^{k} =y^{k+1},\quad \tilde{\lambda }^{k}=\lambda ^{k}-\beta \left ( Ax^{k+1}+By^{k}-b   \right ).    
\end{equation}
As a result, it is feasible to rewrite the inequalities \eqref{eqn:131} and \eqref{eqn:132} respectively as
\begin{subequations}\label{eqn:134}
\begin{align}
     &\tilde{x}^{k} \in X ,\quad \theta _{1} \left ( x \right )  -\theta _{1} \left ( \tilde{x} ^{k}  \right ) +\left ( x-\tilde{x} ^{k}  \right )^{T}\left ( -A^{T}\tilde{\lambda } ^{k}   \right ) \ge 0.\quad \forall x\in X, \label{eqn:134:a}\\
\intertext{and}
\tilde{y } ^{k}\in Y ,&\quad \theta _{2} \left ( y \right ) -\theta _{2} \left ( \tilde{y} ^{k}  \right )+\left ( y-\tilde{y}^{k}   \right )^{T}\left \{ -B^{T}\tilde{\lambda }^{k} +\beta   \delta _{k}\left ( \tilde{y}^{k}-y ^{k}     \right )     \right \}\ge 0,\quad\forall y\in Y.\label{eqn:134:b}
\intertext{Note that $\tilde{\lambda } ^{k}$, defined in \eqref{eqn:133}, may alternatively be portrayed as the variational inequality:}
&\tilde{\lambda }^{k}\in R^{m}, \quad  \left ( \lambda -\tilde{\lambda } ^{k}  \right ) ^{T}\left \{ A\tilde{x} ^{k}+B\tilde{y}^{k} -b - B\left ( \tilde{y} ^{k}-y  ^{k}   \right ) +\frac{ 1 }{\beta } \left ( \tilde{\lambda }^{k}-\lambda  ^{k}    \right )   \right \}\notag \\&\ge 0,\quad\forall \lambda \in R^{l}.\label{eqn:134:c}
\end{align}
\end{subequations}
Using the notation from \eqref{eqn:114a} and the matrix $Q_{k+1}$ described in \eqref{eqn:117}, we can rewrite the inequalities \eqref{eqn:134:a}-\eqref{eqn:134:c} as the variational inequality:
$$\tilde{w} ^{k} \in \Omega ,\quad \theta \left ( u \right ) -\theta \left ( \tilde{u} ^{k}  \right ) +\left ( w-\tilde{w}^{k}   \right )^{T} F\left ( \tilde{w}^{k} \right )  \ge \left ( v-\tilde{v}^{k}   \right )^{T}  Q_{k+1} \left (  v^{k}-\tilde{v}^{k}   \right ),\quad\forall w\in \Omega .$$
Consequently, Lemma 3.1 is established.\\
\textbf{Lemma 3.2}. Let  $ \left \{ w^{k} \right \} $ be the sequence generated by 
 the Algorithm 1 for the problem \eqref{eqn:1} and $ \tilde{w} ^{k}  $ be defined by\eqref{eqn:119} . Then we have 
\begin{equation}\label{eqn:135}
    v^{k} -v^{k+  1} =M  \left ( v^{k}-\tilde{v}^{k}    \right )   
\end{equation}
where $M$ is defined in \eqref{eqn:117}. \\
Proof:  According to the definition of $\tilde{\lambda}^{k}$ in \eqref{eqn:133} and $\lambda ^{k+1} $ in \eqref{eqn:126:3} , we have 
\begin{align*}
\lambda ^{k+  1} &=\lambda ^{k} -\beta \left ( Ax^{k+1}+By^{k+1}-b   \right )\\&= \lambda ^{k} -\beta \left ( Ax^{k+1}+By^{k}-b   \right ) -\beta \left ( By^{k+1}-By^{k}  \right )\\&=\tilde{\lambda }^{k}-\beta \left ( By^{k+1}-By^{k}   \right )\\&=\lambda ^{k}-\left ( \lambda ^{k}-\tilde{\lambda }^{k}    \right )+  \beta B\left ( y^{k}-\tilde{y}^{k}    \right ) .    
\end{align*}
Recall $y^{k+1} =\tilde{y}^{k}$ and the notation in \eqref{eqn:133}, we have 
$$\begin{pmatrix}y^{k+1} 
 \\\lambda ^{k+1} 
\end{pmatrix} = \begin{pmatrix}y^{k} 
 \\\lambda ^{k} 

\end{pmatrix}- \begin{pmatrix}
  I_{n_{2} } &O \\
  -\beta B&I_{m }
\end{pmatrix}  \begin{pmatrix}y^{k}-\tilde{y}^{k}   
 \\\lambda ^{k}-\tilde{\lambda }^{k}   
\end{pmatrix}.$$ \\
Consequently, Lemma 3.2 is established. \\
\textbf{Lemma 3.3}. Let  $ \left \{ w^{k} \right \} $ be the sequence generated by 
 the Algorithm 1 for the problem \eqref{eqn:1} , then for any $w^{\ast } =\left ( x^{\ast },y^{\ast } ,\lambda ^{\ast } \right )\in \Omega ^{\ast },$ we have
 \begin{equation}\label{eqn:136}
\begin{split}
 &\delta _{k}\beta   \left ( y^{k+1}-y^{\ast }  \right )^{T}\left ( y^{k+1} -y^{k}  \right ) -\beta   \left ( y^{k+1}-y^{\ast }   \right )^{T}B^{T} B\left ( y^{k+1} -y^{k}  \right )\\& + \frac{1}{\beta }\left ( \lambda ^{k}-\lambda ^{k+1}   \right )  ^{T} \left ( \lambda ^{\ast } -\lambda ^{k+1}  \right ) \\&\le \beta \left ( Ax^{k+1} -Ax^{\ast } \right )^{T}\left ( By^{k+1}-By^{k}  \right ) .    
\end{split}
\end{equation}
Proof: Since $w^{\ast }$ is a solution of $ VI\left ( \Omega ,F,\theta \right )$ and $x^{k+1}\in X,y^{k+1}\in Y,$   we have:
\begin{equation}\label{eqn:137}
\theta _{1} \left ( x^{k+1}  \right ) -\theta _{1} \left ( x^{\ast }  \right ) +\left ( x^{k+1}-x^{\ast }   \right )^{T} \left ( -A^{T}\lambda ^{\ast }   \right )\ge 0  .    
\end{equation}
\begin{equation}\label{eqn:138}
\theta _{2} \left ( y^{k+1}  \right ) -\theta _{2} \left ( y^{\ast }  \right ) +\left ( y^{k+1}-y^{\ast }   \right )^{T} \left ( -B^{T}\lambda ^{\ast }   \right )\ge 0  .    
\end{equation}  
Setting $x=x^{\ast } $ in \eqref{eqn:131} and for any $ x^{k+1}\in X $, we have 
$$\theta _{1} \left ( x^{\ast }  \right )  -\theta _{1} \left ( x ^{k+1}  \right ) +\left ( x^{\ast }-x ^{k+1}  \right )^{T}\left \{ -A^{T}\lambda  ^{k}  +\beta A^{T}\left ( Ax^{k+1} +By^{k} -b   \right ) \right \} \ge 0.$$
According to the definition of $\lambda ^{k+1} $ in \eqref{eqn:126:3}, the above inequality  may be expressed as:
\begin{equation}\label{eqn:139}
\theta _{1} \left ( x^{\ast } \right )  -\theta _{1} \left ( x ^{k+1}  \right ) +\left ( x^{\ast }-x ^{k+1}  \right )^{T}\left \{ -A^{T}\lambda  ^{k+1} +  \beta A^{T} B\left ( y^{k}-y^{k+1}   \right )  \right \} \ge 0.   
\end{equation}
Setting $y=y^{\ast } $ in \eqref{eqn:132} and for any $ y^{k+1}\in Y $, we have 
$$\theta _{2} \left ( y^{\ast } \right ) -\theta _{2} \left ( y ^{k+1}  \right )+\left ( y^{\ast }-y^{k+1}   \right )^{T}\left \{-B^{T}\lambda ^{k}+  \beta B^{T}\left ( Ax^{k+1}+  By^{k} -b   \right )\right.$$ $$\left.+  \beta  \delta _{k}\left ( y^{k+1}-y^{k}    \right )  \right \}\ge 0.$$
According to the definition of $\lambda ^{k+1} $ in \eqref{eqn:126:3}, the above inequality   may be expressed as:
\begin{equation}\label{eqn:140}
 \begin{split}
\theta _{2} \left ( y^{\ast } \right ) -\theta _{2} \left ( y ^{k+1}  \right )&+\left ( y^{\ast }-y^{k+1}   \right )^{T}\left \{ -B^{T}\lambda ^{k+1} +\beta  \delta _{k}    \left ( y^{k+1}-y ^{k}\right )\right.\\& -\left.\beta B^{T}B   \left ( y^{k+1}-y ^{k}     \right )            \right \}\ge 0.     
 \end{split}   
\end{equation}
Adding \eqref{eqn:137} and \eqref{eqn:139}, we can get 
$$\left ( x^{k+1} -x^{\ast }  \right )^{T}\left \{ A^{T}\lambda ^{k+1} -A^{T}\lambda ^{\ast } -\beta A^{T}B\left ( y^{k}-y^{k+1}   \right )     \right \}\ge 0.$$
The above inequality may be expressed as:
\begin{equation}\label{eqn:141}
\left ( Ax^{k+1} -Ax^{\ast }  \right )^{T}\left ( \lambda ^{\ast }- \lambda ^{k+1}\right )\le \beta \left ( Ax^{k+1} -Ax^{\ast } \right )^{T}\left ( By^{k+1}-By^{k}  \right )  .   
\end{equation}
Adding \eqref{eqn:138} and \eqref{eqn:140}, we can get
\begin{equation}\label{eqn:142}
\begin{split}
 \delta _{k}\beta   \left ( y^{k+1}-y^{\ast }   \right )^{T}\left ( y^{k+1} -y^{k}  \right ) -\beta   \left ( y^{k+1}-y^{\ast }   \right )^{T}B^{T} B\left ( y^{k+1} -y^{k}  \right )  \\+  \left ( By^{k+1} -By^{\ast }  \right ) ^{T} \left ( \lambda ^{\ast } -\lambda ^{k+1}  \right ) \le 0 .    
\end{split}
\end{equation}
Adding \eqref{eqn:141}and \eqref{eqn:142} , we can get
\begin{equation}\label{eqn:143}
\begin{split}
& \delta _{k}\beta   \left ( y^{k+1}-y^{\ast }   \right )^{T}\left ( y^{k+1} -y^{k}  \right ) -\beta   \left ( y^{k+1}-y^{\ast }   \right )^{T}B^{T} B\left ( y^{k+1} -y^{k}  \right )\\& + \left (Ax^{k+1} + By^{k+1} -b  \right ) ^{T} \left ( \lambda ^{\ast } -\lambda ^{k+1}  \right ) \\&\le \beta \left ( Ax^{k+1} -Ax^{\ast } \right )^{T}\left ( By^{k+1}-By^{k}  \right )  .     
\end{split}
\end{equation}
For the dual subproblem \eqref{eqn:126:3}, we have 
\begin{equation}\label{eqn:144}
 Ax^{k+1} + By^{k+1} -b=\frac{1}{\beta}\left ( \lambda ^{k}-\lambda ^{k+1} \right ) ,  
\end{equation}
Substituting \eqref{eqn:144} into \eqref{eqn:143} , we get 
\begin{equation*}
\begin{split}
 &\delta _{k}\beta   \left ( y^{k+1}-y^{\ast }  \right )^{T}\left ( y^{k+1} -y^{k}  \right ) -\beta   \left ( y^{k+1}-y^{\ast }   \right )^{T}B^{T} B\left ( y^{k+1} -y^{k}  \right )\\& + \frac{1}{\beta }\left ( \lambda ^{k}-\lambda ^{k+1}   \right )  ^{T} \left ( \lambda ^{\ast } -\lambda ^{k+1}  \right ) \\&\le \beta \left ( Ax^{k+1} -Ax^{\ast } \right )^{T}\left ( By^{k+1}-By^{k}  \right ) .    
\end{split}
\end{equation*}
Consequently, Lemma 3.3 is established.\\
\textbf{Lemma 3.4}.  Let  $ \left \{ w^{k} \right \} $ be the sequence generated by 
 the Algorithm 1 for the problem \eqref{eqn:1} , then for any $w^{\ast } =\left ( x^{\ast },y^{\ast } ,\lambda ^{\ast } \right )\in \Omega ^{\ast },$ we have
\begin{equation}\label{eqn:145}
\begin{split}
&\delta _{k}\beta \left \| y^{k+1}-y^{\ast }    \right \|^{2} +  \frac{1}{\beta }\left \| \lambda ^{k+1}-\lambda ^{\ast }    \right \|^{2}   \\& \le\delta _{k}\beta \left \| y^{k}-y^{\ast }    \right \|^{2} +  \frac{1}{\beta }\left \| \lambda ^{k}-\lambda ^{\ast }    \right \|^{2} \\&\quad-
\left ( \delta _{k}\beta \left \| y^{k+1}-y^{k }    \right \|^{2} +  \frac{1}{\beta }\left \| \lambda ^{k+1}-\lambda ^{k }    \right \|^{2}  \right ) \\&\quad-2 \left ( \lambda ^{k+1}-\lambda ^{k}  \right ) ^{T} \left ( By^{k+1} - By^{k}\right ).
\end{split}    
\end{equation}
Proof: Taking advantage of the identity:$$\left \| a+b \right \|^{2}= \left \| a \right \|^{2} -\left \| b \right \|^{2}+  2\left ( a+b \right )^{T}b .$$ 
We are able to get:
\begin{equation}\label{eqn:146}
\begin{split}
&\delta _{k}\beta \left \| y^{k+1}-y^{\ast }    \right \|^{2} +  \frac{1}{\beta }\left \| \lambda ^{k+1}-\lambda ^{\ast }    \right \|^{2}+  \beta \left \| B\left ( y^{k+1}-y^{\ast }   \right )  \right \|^{2} \\& =\delta _{k}\beta \left \| y^{k}-y^{\ast }    \right \|^{2} +  \frac{1}{\beta }\left \| \lambda ^{k}-\lambda ^{\ast }    \right \|^{2}+  \beta \left \| B\left ( y^{k}-y^{\ast }   \right )  \right \|^{2}\\&\quad-
\left (  \delta _{k}\beta \left \| y^{k+1}-y^{k }    \right \|^{2} +  \frac{1}{\beta }\left \| \lambda ^{k+1}-\lambda ^{k }    \right \|^{2}+  \beta \left \| B\left ( y^{k+1}-y^{k }   \right )  \right \|^{2}\right ) \\& \quad +2 \delta _{k}\beta   \left ( y^{k+1}-y^{\ast }  \right )^{T}\left ( y^{k+1} -y^{k}  \right ) + \frac{2}{\beta }\left ( \lambda ^{k}-\lambda ^{k+1}   \right )  ^{T} \left ( \lambda ^{\ast } -\lambda ^{k+1}  \right )\\&\quad+2\beta \left ( By^{k+1}-By^{\ast }   \right )^{T}\left ( By^{k+1} -By^{k} \right ) .      
\end{split}    
\end{equation}
\begin{equation}\label{eqn:147}
\begin{split}
\quad-\beta \left \| B\left ( y^{k+1}-y^{\ast }   \right )  \right \|^{2}=&- \beta \left \| B\left ( y^{k}-y^{\ast }   \right )  \right \|^{2}+ \beta \left \| B\left ( y^{k+1}-y^{k }   \right )  \right \|^{2}\\&-2\beta   \left ( y^{k+1}-y^{\ast }  \right )^{T}B^{T} B\left ( y^{k+1} -y^{k}  \right ).
\end{split}    
\end{equation}
Adding \eqref{eqn:146} and \eqref{eqn:147}, it indicates that
\begin{equation}\label{eqn:148}
\begin{split}
&\delta _{k}\beta \left \| y^{k+1}-y^{\ast }    \right \|^{2} +  \frac{1}{\beta }\left \| \lambda ^{k+1}-\lambda ^{\ast }    \right \|^{2}   \\& =\delta _{k}\beta \left \| y^{k}-y^{\ast }    \right \|^{2} +  \frac{1}{\beta }\left \| \lambda ^{k}-\lambda ^{\ast }    \right \|^{2} \\&\quad-
\left (  \delta _{k}\beta \left \| y^{k+1}-y^{k }    \right \|^{2} +  \frac{1}{\beta }\left \| \lambda ^{k+1}-\lambda ^{k }    \right \|^{2}  \right ) \\& \quad +2  \delta _{k}\beta   \left ( y^{k+1}-y^{\ast }  \right )^{T}\left ( y^{k+1} -y^{k}  \right )-2\beta   \left ( y^{k+1}-y^{\ast }  \right )^{T}B^{T} B\left ( y^{k+1} -y^{k}  \right ) \\&\quad+ \frac{2}{\beta }\left ( \lambda ^{k}-\lambda ^{k+1}   \right )  ^{T} \left ( \lambda ^{\ast } -\lambda ^{k+1}  \right )+2\beta \left ( By^{k+1}-By^{\ast }   \right )^{T}\left ( By^{k+1} -By^{k} \right ) .      
\end{split}    
\end{equation}
Next, we consider the crossing terms in \eqref{eqn:148} . \\
If we combine $Ax^{\ast } +By^{\ast } = b$ with \eqref{eqn:136}, we obtain
\begin{equation}\label{eqn:149}
\begin{split}
&\quad2 \delta _{k}\beta   \left ( y^{k+1}-y^{\ast }  \right )^{T}\left ( y^{k+1} -y^{k}  \right )-2\beta   \left ( y^{k+1}-y^{\ast }  \right )^{T}B^{T} B\left ( y^{k+1} -y^{k}  \right ) \\&+ \frac{2}{\beta }\left ( \lambda ^{k}-\lambda ^{k+1}   \right )  ^{T} \left ( \lambda ^{\ast } -\lambda ^{k+1}  \right )+2\beta \left ( By^{k+1}-By^{\ast }   \right )^{T}\left ( By^{k+1} -By^{k} \right ) \\&\le 2\beta \left ( Ax^{k+1} -Ax^{\ast } \right )^{T}\left ( By^{k+1}-By^{k}  \right )\\&\quad+2\beta \left ( By^{k+1}-By^{\ast }   \right )^{T}\left ( By^{k+1} -By^{k} \right ) 
\\&=2\beta \left ( Ax^{k+1}+By^{k+1} -b  \right ) ^{T}\left ( By^{k+1}-By^{k}  \right )\\& =-2 \left ( \lambda ^{k+1}-\lambda ^{k}  \right ) ^{T} \left ( By^{k+1} - By^{k}\right ).
\end{split}    
\end{equation}
Substituting \eqref{eqn:149} into \eqref{eqn:148}, the assertion of this lemma is proved directly.\\
\textbf{Theorem 3.1}. Let  $ \left \{ w^{k} \right \} $ be the sequence generated by 
 the Algorithm 1 for the problem \eqref{eqn:1} and $w^{\ast } =\left ( x^{\ast },y^{\ast } ,\lambda ^{\ast } \right )\in \Omega ^{\ast }$ is a solution of $VI\left ( \Omega ,F,\theta \right )$,  then for $\forall \varepsilon \in \left ( 0, \frac{1}{2} \right )$ and $D_{k+1} = \delta _{k}I_{n_{2} } - \frac{1}{2\varepsilon }B^{T}B\succ 0 $, we have
\begin{equation}\label{eqn:150}
\begin{split}
\left \| v^{k+1} -v^{\ast }  \right \|_{H_{k+1} }^{2}&\le\left \| v^{k} -v^{\ast }  \right \|_{H_{k+1} }^{2}\\&- 
\left ( \beta \left \| y^{k+1}-y^{k }    \right \| _{D_{k+1} }^{2} +  \frac{1-2\varepsilon }{\beta }\left \| \lambda ^{k+1}-\lambda ^{k }    \right \|^{2}  \right ) .  
\end{split}    
\end{equation}
Proof: According to  the $Cauchy-Schwarz$ inequality, for $\forall  \varepsilon \in \left ( 0, \frac{1}{2} \right )$, we have 
\begin{equation}\label{eqn:151}
\left ( \lambda ^{k+1}-\lambda ^{k}  \right ) ^{T} \left ( By^{k+1} - By^{k}\right )\ge -\frac{\varepsilon }{\beta }\left \| \lambda ^{k+1}- \lambda ^{k}  \right \|^{2}-\frac{1}{4\varepsilon }\beta \left \| B\left ( y^{k+1}-y^{k}   \right )  \right \|^{2} .
\end{equation}
Substituting \eqref{eqn:151} into \eqref{eqn:145} , we get 
\begin{equation*}
\begin{split}
&\delta _{k}\beta \left \| y^{k+1}-y^{\ast }    \right \|^{2} +  \frac{1}{\beta }\left \| \lambda ^{k+1}-\lambda ^{\ast }    \right \|^{2}   \\& \le \delta _{k}\beta \left \| y^{k}-y^{\ast }    \right \|^{2} +  \frac{1}{\beta }\left \| \lambda ^{k}-\lambda ^{\ast }    \right \|^{2} \\&\quad-
\left (  \delta _{k}\beta \left \| y^{k+1}-y^{k }    \right \|^{2} +  \frac{1}{\beta }\left \| \lambda ^{k+1}-\lambda ^{k }    \right \|^{2}  \right ) \\&\quad+\frac{2\varepsilon }{\beta }\left \| \lambda ^{k+1}- \lambda ^{k}  \right \|^{2}+\frac{1}{2\varepsilon }\beta \left \| B\left ( y^{k+1}-y^{k}   \right )  \right \|^{2}
\end{split}    
\end{equation*}
The above inequality is equivalent to:
\begin{equation}\label{eqn:152}
\begin{split}
&\delta _{k}\beta \left \| y^{k+1}-y^{\ast }    \right \|^{2} +  \frac{1}{\beta }\left \| \lambda ^{k+1}-\lambda ^{\ast }    \right \|^{2}   \\& \le\delta _{k}\beta \left \| y^{k}-y^{\ast }    \right \|^{2} +  \frac{1}{\beta }\left \| \lambda ^{k}-\lambda ^{\ast }    \right \|^{2} \\&\quad-
\left ( \beta \left \| y^{k+1}-y^{k }    \right \| _{D_{k+1} }^{2} +  \frac{1-2\varepsilon }{\beta }\left \| \lambda ^{k+1}-\lambda ^{k }    \right \|^{2}  \right ) .
\end{split}    
\end{equation}
 where $D_{k+1} \succ 0 $ is equivalent to $\left \| y^{k+1}-y^{k }    \right \| _{D_{k+1} }^{2}> 0.$\\
 then it is only necessary to $\delta _{k}\beta\left \| y^{k+1}-y^{k }    \right \|^{2}-\frac{1}{2\varepsilon }\beta \left \| B\left ( y^{k+1}-y^{k}   \right )  \right \|^{2}>0.$\\
 According to the definition of $H_{k+1} $ in \eqref{eqn:118}, inequality \eqref{eqn:152} is equivalent to:
\begin{equation*}
\begin{split}
\left \| v^{k+1} -v^{\ast }  \right \|_{H_{k+1} }^{2}&\le\left \| v^{k} -v^{\ast }  \right \|_{H_{k+1} }^{2}\\&- 
\left ( \beta \left \| y^{k+1}-y^{k }    \right \| _{D_{k+1} }^{2} +  \frac{1-2\varepsilon }{\beta }\left \| \lambda ^{k+1}-\lambda ^{k }    \right \|^{2}  \right ) .  
\end{split}    
\end{equation*}
Therefore, Theorem 3.1 is proved .\\
\textbf{Theorem 3.2}.  Let  $ \left \{ w^{k} \right \} $ be the sequence generated by 
 the Algorithm 1 for the problem \eqref{eqn:1} , then for any $w^{\ast } =\left ( x^{\ast },y^{\ast } ,\lambda ^{\ast } \right )\in \Omega ^{\ast },$ we have
\begin{equation}\label{eqn:153}
\begin{split}
 \left \| v^{k+1} -v^{\ast }  \right \|_{H_{k+1} }^{2}&\le \left ( 1+  \xi _{k}\right ) \left \| v^{k} -v^{\ast }  \right \|_{H_{k} }^{2}\\&\quad- 
\left ( \beta \left \| y^{k+1}-y^{k }    \right \| _{D_{k+1} }^{2} +  \frac{1-2\varepsilon }{\beta }\left \| \lambda ^{k+1}-\lambda ^{k }    \right \|^{2}  \right ).
\end{split}
\end{equation}
where $\xi _{k}=\max\left \{ 0,\delta _{k}-\delta _{k-1}   \right \} .$\\
Proof: According to \eqref{eqn:150} and the definition of $H_{k+1}$ in \eqref{eqn:118}, we have
\begin{equation}\label{eqn:154}
\begin{split}
 \left \| v^{k+1} -v^{\ast }  \right \|_{H_{k+1} }^{2}&\le\left \| v^{k} -v^{\ast }  \right \|_{H_{k+1} }^{2}\\&- 
\left ( \beta \left \| y^{k+1}-y^{k }    \right \| _{D_{k+1} }^{2} +  \frac{1-2\varepsilon }{\beta }\left \| \lambda ^{k+1}-\lambda ^{k }    \right \|^{2}  \right ).    
\end{split}
\end{equation}
 Due to
\begin{equation}\label{eqn:155}
\begin{split}
\left \| v^{k} -v^{\ast }  \right \|_{H_{k+1} }^{2}
&= \left \| v^{k} -v^{\ast }  \right \|_{H_{k} }^{2}+\left ( \delta _{k}-\delta _{k-1}  \right )\left \| y^{k} -y^{\ast }  \right \|^{2} \\&\le \left ( 1+  \xi _{k}   \right ) \left \| v^{k} -v^{\ast }  \right \|_{H_{k} }^{2}.
\end{split}
\end{equation}
Substituting \eqref{eqn:155} into \eqref{eqn:154} , we can get
\begin{equation*}
\begin{split}
 \left \| v^{k+1} -v^{\ast }  \right \|_{H_{k+1} }^{2}&\le \left ( 1+  \xi _{k} \right ) \left \| v^{k} -v^{\ast }  \right \|_{H_{k} }^{2}\\&\quad- 
\left ( \beta \left \| y^{k+1}-y^{k }    \right \| _{D_{k+1} }^{2} +  \frac{1-2\varepsilon }{\beta }\left \| \lambda ^{k+1}-\lambda ^{k }    \right \|^{2}  \right ).
\end{split}
\end{equation*}
Therefore,  Theorem 3.2 is proved.\\
According to the definition of $\delta _{k}$ in Algorithm 1, we can get 
$$\sum_{k=0}^{+  \infty } \xi_{k} \le \left \lceil log_{\eta }\left ( \frac{\left \| B^{T}B  \right \| }{\delta _{min} }  \right )   \right \rceil\left \| B^{T}B  \right \| < +  \infty .$$
where $\left \lceil x \right \rceil$ is the smallest integer greater than or equal to $x$ for any  $x\in R $.\\
\textbf{Theorem 3.3}.  Let  $ \left \{ w^{k} \right \} $ be the sequence generated by 
 the Algorithm 1 for the problem \eqref{eqn:1} , then for any $w^{\ast } =\left ( x^{\ast },y^{\ast } ,\lambda ^{\ast } \right )\in \Omega ^{\ast },$ we have\\
(a).$ \lim_{ k\to \infty}\left \| v^{k+1} -v^{\ast }  \right \|_{H_{k+1} }^{2}$ exists and $\lim_{ k\to \infty}  \left \| v^{k+1} -v^{\ast }  \right \|_{H_{k+1} }^{2} < +  \infty.$ \\ 
(b).$\lim_{ k\to \infty}\left \| y^{k+1}-y^{k}   \right \|= 0$ and $\lim_{ k\to \infty}\left \| \lambda ^{k+1}-\lambda ^{k}   \right \|= 0 .$  \\
Proof: Letting $$a^{k+1}= \left \| v^{k+1} -v^{\ast }  \right \|_{H_{k+1} }^{2},\quad b^{k} = \xi _{k},\quad c^{k} =0 ,$$$$d^{k}=\beta \left \| y^{k+1}-y^{k }    \right \| _{D_{k+1} }^{2} +  \frac{1-2\varepsilon }{\beta }\left \| \lambda ^{k+1}-\lambda ^{k }    \right \|^{2}.$$\\
According to Lemma 2.2, Theorem 3.3 is proved. \\
\textbf{Theorem 3.4}.  Let  $ \left \{ w^{k} \right \} $ be the sequence generated by 
 the Algorithm 1 for the problem \eqref{eqn:1} , then $\left \{ w^{k}  \right \} $ converges to  $w^{\infty } \in
\Omega ^{\ast }$.\\
Proof: According to Theorem 3.3, we can get 
\begin{equation}\label{eqn:156}
\lim_{k \to \infty}\left \| v^{k+1}-v^{k}   \right \|=0.    
\end{equation}
 Combining \eqref{eqn:135} and the non-singularity of $M$, we have
 \begin{equation}
  \lim_{k \to \infty} \left \| v^{k}-\tilde{v}^{k}    \right \|=0 .   
 \end{equation}
According to Theorem 3.3,  the sequence $\left \{ \left \| v^{k+1} -v^{\ast }  \right \|_{H_{k+1} }^{2} \right \}$ is a bounded sequence, \\ then for any fixed $v^{\ast }\in \Omega ^{\ast } $, $\left \| v^{k+1}-v^{\ast }   \right \|$  is bounded.\\
which means that the sequence $\left \{ v^{k} \right \}$ is also bounded. \\
Because of $$\left \| \tilde{v}^{k}-v^{\ast }    \right \|\le \left \| v^{k}-\tilde{v}^{k}    \right \|+\left \| v^{k}-v^{\ast }   \right \|.$$
It is known that $\left \| \tilde{v}^{k}-v^{\ast }    \right \|$ is bounded. \\
Clearly the sequence $\left \{ \tilde{v}^{k}   \right \} $ is also bounded and there must exist a convergence point $v^{\infty }, $ such that a subsequence $\left \{ \tilde{v}^{k_{j} } \right \}$ of the existence sequence $\left \{ \tilde{v}^{k} \right \}$  converges to $v^{\infty } . $\\
According to  $\tilde{\lambda }^{k}=\lambda ^{k} -\beta \left ( A\tilde{x}^{k}+  By^{k}-b    \right ) $ in \eqref{eqn:133}, we have
\begin{equation}\label{eqn:158}
A\tilde{x}^{k_{j} } =\frac{1}{\beta } \left ( \lambda ^{k_{j} } -\tilde{\lambda }^{k_{j} }   \right ) -\left ( By^{k_{j} }-b  \right ).    
\end{equation}
Since Matrix A is a column full rank matrix, we can get that the sequence $\left \{ \tilde{x}^{k_{j} }   \right \}$ converges. \\Set $ \lim_{k_{j}  \to \infty} \tilde{x}^{k_{j} } =x^{\infty } ,$ then there exists a subsequence $\left \{ \tilde{w}^{k_{j} }   \right \} $ converging to $w^{\infty }$.  \\  
Setting $ k=k_{j} $ in \eqref{eqn:130}, we have
$\tilde{w} ^{k_{j} }\in \Omega $ such that\\$$\theta \left ( u \right )-  \theta \left ( \tilde{u}^{k_{j} }  \right )+  \left ( w-\tilde{w}^{k_{j} }   \right )^{T}F\left ( \tilde{w}^{k_{j} }  \right )\ge \left ( w-\tilde{w}^{k_{j} }  \right ) ^{T}Q\left ( w^{k_{j} } -\tilde{w}^{k_{j} }  \right ) ,\quad\forall w\in \Omega . $$\\
Let the above inequality $k \to \infty ,$ it is clear that:
\begin{equation}\label{eqn:159}
 w^{\infty } \in \Omega ,\quad\theta \left ( u \right ) -\theta \left ( u^{\infty }  \right ) +  \left ( w-w^{\infty }  \right )^{T}F\left ( w^{\infty }  \right )\ge 0,\quad\forall w\in \Omega .    
\end{equation}
 \eqref{eqn:159} shows that $w^{\infty } $ in $\Omega ^{\ast }$ is a solution of the variational inequality $VI\left ( \Omega ,F,\theta \right )$, so the sequence $\left \{ w^{k} \right \} $ converges to $w^{\infty }. $

\section{Numerical Experiments}\label{sec6}
\noindent This work presents the numerical performance of Algorithm 1 for solving the Lasso problem. We compare Algorithm 1 to the OLADMM presented by \cite{28}. All simulations are run on a laptop with 4GB RAM as well as Anaconda 3.\\
We initially study the least absolute shrinkage and selection operator (LASSO) problem. (see,e.g.,\cite{31}):
\begin{equation}\label{eqn:160}
    \min\left \{ \frac{1}{2} \left \| Ay-b \right \|^{2} +\sigma \left \| y \right \|_{1}|y\in R^{n}     \right \} ,
\end{equation}
where $A\in R^{m\times n}$  with $ m\ll n$, $\left \| y \right \| _{1} =\sum_{i=1}^{n}\left | y_{i}  \right |$ , $b\in R^{m}$ and $\sigma > 0$ is a regularization parameter. This model may be interpreted as choosing main features among $n$ give
ones and there are $m$ data points; and $b$ is understood as labels.

To solve the model \eqref{eqn:160} by various ADMM-based schemes, one way is introducing an auxiliary variable $x$ and the constraint $x=Ay$, and then rewritting \eqref{eqn:160} into the form of \eqref{eqn:1}:
\begin{equation}\label{eqn:161}
   \min_{x,y}\left \{  \frac{1}{2}\left \| x-b \right \|_{2}^{2}   + \sigma     \left \| y \right \|_{1} |  x=Ay, x\in R^{m} , y\in R^{n}  \right \}  
\end{equation}
Then, the implementation of the Algorithm 1 requires computing the variables $x$ and $y$ , respectively, by the following two schemes:
\begin{equation}
x^{k+1} =\frac{1}{1+  \beta }\left ( b+  \lambda ^{k}+  \beta Ay^{k}   \right )  
\end{equation}
 and
 \begin{equation}
   y^{k+1} =shrink\left \{ y^{k}-\frac{1}{\delta_{k}\beta   }A^{T} \left [ \lambda ^{k}-\beta \left ( x^{k+1}-Ay^{k}   \right )   \right ]  ,\frac{\sigma }{ \delta_{k}\beta   }    \right \}  
 \end{equation}
We use the same data creation technique as in \cite{31}: first choose $A_{i,j} \sim N\left (0,1 \right)$, then normalize the columns to have unit $l_{2} $ norm. The true value $y_{true}$ is constructed using a density of $100/n$ non-zero entries sampled from a $N(0,1)$ distribution. We create the labels $b$ as $b=Ay_{true} +\varrho $, where $\varrho $ is a random noise with $\varrho \sim N\left (0,10^{-3}I \right) $, and we set the regularization parameter $\sigma =0.1\sigma _{max}$ with $ \sigma _{max} =\left \| A^{T} b \right \| _{\infty } $.
Define each algorithm's parameters as follows:\\
Algorithm 1: $\delta_{0} =\delta_{-1}=0.75\left \| A^{T}A \right \|$,$\delta _{min}=0.05\left \| B^{T}B  \right \|  $, $\tau   =1.1$  ,$\eta =1.1$, $\beta =1$, $\varepsilon =\frac{5}{11} $
.\\
$\text{OLADMM: } r=0.75\left \| A^{T}A  \right \|,  \beta =1.$\\
The stopping criterion is: 
\begin{equation}
 \left \| p^{k}  \right \| =\left \| x^{k+1} -Ay^{k+1}  \right \| < \epsilon  ^{pri} , \quad and \quad \left \| q^{k}  \right \|= \left \| \beta A \left ( y^{k+1}-y^{k}   \right )   \right \| < \epsilon  ^{dual} ,   
\end{equation}
where 
$$\epsilon ^{pri}=\sqrt{n}\epsilon  ^{abs} +\epsilon  ^{rel} max\left \{ \left \| x^{k+1}  \right \| ,\left \| Ay^{k+1}  \right \|  \right \} ,\quad \epsilon  ^{dual} =\sqrt{n}\epsilon  ^{abs} +\epsilon  ^{rel} \left \| y^{k+1}  \right \|.$$with $\epsilon  ^{abs}$ and $\epsilon  ^{rel}$ as $10^{-6} $ and $10^{-4} .$ The initial points $\left ( y^{0},\lambda ^{0}   \right )$ ia set as zero. \\
\text{Table 1.} Comparsion between Algorithm 1 and OLADMM for \eqref{eqn:161}.  \\
\begin{tabular}{|c|c|c|c|c|c|c|c|c|c|}  
\hline \multicolumn{2}{|c|}{$n \times n$ matrix } & \multicolumn{4}{|c|}{ OLADMM } & \multicolumn{4}{|c|}{Algorithm 1} \\  
\hline$m$ & $n$ & Iter. & CPU(s) & $\left\|p^{k}\right\|$ & $\left\|q^{k}\right\|$ & Iter. & CPU(s) & $\left\|p^{k}\right\|$ & $\left\|q^{k}\right\|$ \\    
\hline 1000 & 1500 & 404 & 9.93 & $0.00155 $ & $0.00155$ & 47 & 1.62 & $0.00137$ & $0.00137$ \\   
\hline 1000 & 2000 & 456 & 14.80 & $0.00161$ & $0.00161$ & 50& 2.61 & $0.00157$ & $0.00157$ \\
\hline 1500 & 3000 & 484 & 33.49 & $0.00159$ & $0.00159$ & 55 & 5.70 & $0.00145$ & $0.00145$ \\  
\hline 2000 & 3000 & 422 & 39.36 & $0.00163$ & $0.00163$ & 45 & 6.61 & $0.00150$ & $0.00150$ \\  
\hline 2000 & 4000 & 479 & 60.12 & $0.00154$ & $0.00154$ & 51 & 10.35 & $0.00153$ & $0.00153$ \\ 
\hline 3000 & 4000 & 403 & 76.78 & $0.00157$ & $0.00157$ & 43& 14.18 & $0.00133$ & $0.00133$ \\  
\hline 3000 & 5000 & 455 & 107.78 & $0.00146$ & $0.00146$ & 50& 20.51 & $0.00128$ & $0.00128$ \\
\hline 4000 & 5000 & 415 & 131.88 & $0.00147$ & $0.00147$ & 45& 26.41 & $0.00128$ & $0.00128$ \\
\hline
\end{tabular}\\

Table 1 shows the number of iterations and overall calculation time required by OLADMM and Algorithm 1. The table 1 shows that Algorithm 1 is quicker than OLADMM, requiring fewer iteration steps and taking less time to attain termination tolerance. Consequently, Algorithm 1 is far more effective than OLADMM.  Figures 1 and 2 show the numerical results of two methods for m = 4000 and n = 5500.  Algorithm 1 converges linearly and performs quicker than OLADMM, supporting our convergence findings. Thus, Algorithm 1 superior to OLADMM. Furthermore, Algorithm 1 outperforms OLADMM, providing significant support for our convergence research.

\begin{figure}[htbp]
\centering  
\includegraphics[width=6cm,height = 6cm]{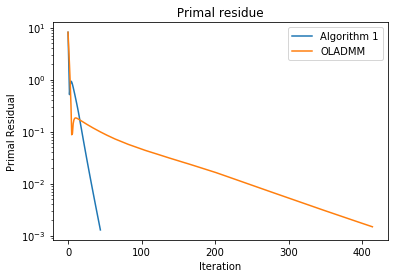}
\includegraphics[width=6cm,height = 6cm]{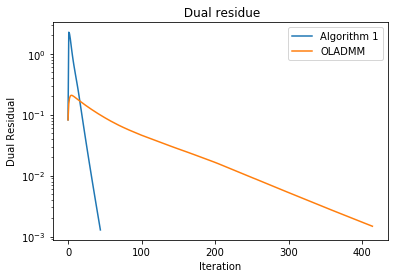}
\caption{ Convergence curves of the primal and dual residuals 
for Algorithm 1 and OLADMM with the number of iterations.}
\end{figure}

\begin{figure}[htbp]
\centering  
\includegraphics[width=6cm,height = 6cm]{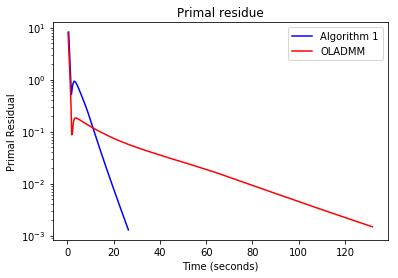}
\includegraphics[width=6cm,height = 6cm]{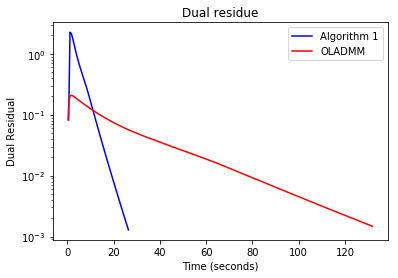}
\caption{ Convergence curves of the primal and dual residuals 
for Algorithm 1 and OLADMM with the time of iterations.}
\end{figure}

\section{Conclusion}\label{secA1}
In this study, we combine the linearized alternating direction multiplier method with an adaptive technique to address the challenge of hard-to-solve subproblems and to improve the efficiency of numerical experiments. We propose an adaptive linearized alternating direction multiplier method and provide an in-depth analysis of the convergence of the algorithm. The results of numerical experiments fully verify the effectiveness of the algorithm, which provides strong methodological and theoretical support for solving complex optimization problems in practical engineering and scientific computing. We believe that this research result has important reference value for academic research and engineering applications in related fields.\\

\textbf{ Conflict of Interest}:The authors declare that they have no conflict of interest.






\bibliography{sn-bibliography}

\end{document}